\documentclass[fleqn]{wiley-article}

\usepackage{siunitx}
\usepackage{supertabular}
\papertype{Original Article}
\paperfield{Journal Section}

\title{Decomposing Convexified Security-Constrained ACOPF Problem with AGC Reformulation}

\abbrevs{}

\begin{document}
\author[\authfn{1}]{Muhammad Waseem}
\author[1\authfn{1}]{Saeed D. Manshadi}


\affil[1]{Department of Electrical \& Computer Engineering, San Diego State University, San Diego, CA 92182 USA}

\corraddress{Saeed D. Manshadi, Department of Electrical \& Computer Engineering, San Diego State University, San Diego, CA 92182 USA}
\corremail{smanshadi@sdsu.edu}



\runningauthor{Muhammad et al.}


\begin{frontmatter}
\maketitle

\begin{abstract}
This paper presents a  reformulation for the  automatic generation control (AGC) in a decomposed convex relaxation algorithm. It finds an optimal solution to the AC optimal power flow (ACOPF) problem that is secure against a large set of contingencies. 
The original ACOPF problem which represents the system without contingency constraints, is convexified by applying the second-order cone relaxation method. 
The contingencies are filtered to distinguish those that will be treated with preventive actions from those that will be left for corrective actions. The  selected contingencies for preventive action are included in the set of security constraints. Benders decomposition is employed to decompose the convexified Security-Constrained ACOPF problem into a master problem and several security check sub-problems. Sub-problems are evaluated in a parallel computing process with enhanced computational efficiency. AGC within each sub-problem is modeled by a set of proposed valid constraints, so the procured solution is the physical response of each generation unit during a contingency. 
Benders optimality cuts are generated for the sub-problems having mismatches and the cuts are passed to the master problem to encounter the security-constraints. 
The accuracy of the relaxation results is verified using the presented tightness measure. The effectiveness of the presented valid AGC constraints and scalability of the proposed algorithm is demonstrated in several case studies.  

\end{abstract}
\end{frontmatter}

\section{Nomenclature}
\vspace{-0.2cm}
\subsection{Indices}
\vspace{0.1cm}
\begin{supertabular}{ll} 
\textit{e} $\in$ $\mathcal{E}$ & Index of line in the set of lines\\
\textit{f} $\in$ $\mathcal{F}$ & Index of transformer in the set of transformers\\
\textit{g} $\in$ $\mathcal{G}$ & Index of generator in the set of generators\\
\textit{k} $\in$ $\mathcal{K}$ & Index of contingency in the set of contingencies\\
\textit{i} $\in$ $\mathcal{I}$ & Index of bus in the set of buses\\
\textit{$i_g$} & Index of bus connected to generator\\
\end{supertabular}
\vspace{-0.4cm}

\subsection{Variables}
\vspace{0.1cm}
\begin{supertabular}{ l p{5in} } 
\textit{\(c_g\)} & Generation cost of generator ($g$) as a function of active power (USD/h)\\
\textit{\(p_e^{ij}\)} & Line \(e\) real power from origin bus \(i\) to bus \(j\)\\
\textit{\(p_f^{ij}\)} & Transformer \(f\) real power from origin bus \(i\) to bus \(j\)\\
\textit{\(p_g\)} & Generator \(g\) real output power\\
\textit{\(q_e^{ij}\)} & Line \(e\) reactive power from origin bus \(i\) to bus \(j\)\\
\textit{\(q_f^{ij}\)} & Transformer \(f\) reactive power from origin bus \(i\) to bus \(j\)\\
\textit{\(q_g\)} & Generator \(g\) reactive output power\\
\textit{\(u_i\)} & Bus \(i\) voltage magnitude\\
\textit{\(e_{i}\)} & Real part of voltage at bus \(i\)\\
\textit{\(f_{j}\)} & Imaginary part of voltage at bus \(j\)\\
\textit{\(\Delta_k\)} & Contingency \(k\) scale factor on generator participation factors defining generator real power contingency response\\
\textit{\(\sigma_{ik}^{P+/-}\)} & Apparent power mismatch for contingency \(k\) in bus \(i\) real power positive/negative parts\\
\textit{\(\sigma_{ik}^{Q+/-}\)} & Apparent power mismatch for contingency \(k\) in bus \(i\) reactive power positive/negative parts\\
\textit{\(\sigma_{ek}\)} & Apparent power mismatch for contingency \(k\) in line \(e\)\\
\textit{\(\sigma_{fk}\)} & Apparent power mismatch for contingency \(k\) in transformer \(f\)\\
\end{supertabular}
\vspace{-0.7cm}

\subsection{Parameters}
\vspace{0.3cm}
\begin{supertabular}{ l  p{5in} }
\textit{\(b_e\)} & Line \(e\) series susceptance\\
\textit{\(b_f\)} & Transformer \(f\) susceptance\\
\textit{\(b_e^{CH}\)} & Line \(e\) total charging susceptance\\
\textit{\(g_e\)} & Line \(e\) series conductance\\
\textit{\(\overline{p}_g\)} & Generator \(g\) real power maximum\\
\textit{\(\underline{p}_g\)} & Generator \(g\) real power minimum\\
\textit{\(\overline{q}_g\)} & Generator \(g\) reactive power maximum\\
\textit{\(\underline{q}_g\)} & Generator \(g\) reactive power minimum\\
\textit{\(\overline{R}_e\)} & Line \(e\) apparent current maximum in base case\\
\textit{\(\overline{R}_e^K\)} & Line \(e\) apparent current maximum in contingencies\\
\textit{\(\overline{s}_f\)} & Transformer \(f\) apparent power maximum in base case\\
\textit{\(\overline{s}_f^K\)} & Transformer \(f\) apparent power maximum in contingencies\\
\textit{\(\overline{u}_i\)} & Bus \(i\) maximum  voltage magnitude in the base case\\
\textit{\(\underline{u}_i\)} & Bus \(i\) minimum voltage magnitude in the base case\\
\textit{\(\overline{u}_i^K\)} & Bus \(i\) maximum voltage magnitude in contingencies\\
\textit{\(\underline{u}_i^K\)} & Bus \(i\) minimum voltage magnitude in contingencies\\
\textit{\(\alpha_g\)} & Participation factor of generator \(g\) in real power contingency response\\
\textit{\(\delta\)} & Weight on base case in objective\\
\textit{\(b_f^M\)} & Transformer \(f\) magnetizing susceptance\\
\textit{\(b_i^{FS}\)} & Bus \(i\) fixed shunt susceptance\\
\textit{\(g_f\)} & Transformer \(f\) series conductance\\
\textit{\(g_f^M\)} & Transformer \(f\) magnetizing conductance\\
\textit{\(g_i^{FS}\)} & Bus \(i\) fixed shunt conductance\\
\textit{\(p_i^L\)} & Bus \(i\) constant real power load\\
\textit{\(q_i^L\)} & Bus \(i\) constant reactive power load\\
\textit{\(\tau_f\)} & Transformer \(f\) tap ratio\\
\textit{\(t_m\)} & Magnitude of complex transformer ratio\\
\textit{\(t_r\)} & Real part of complex transformer ratio\\
\textit{\(t_i\)} & Imaginary part of complex transformer ratio\\
\end{supertabular}
\vspace{-0.4cm}


\section{Introduction}
The inclination to optimally operate the power system may lead to operating points close to the boundary of the security limits. Severe contingencies occurrence in a power system may lead to instabilities or blackouts. Thus, it is  necessary to secure the power system using preventive or corrective control actions. The security assessment is one of the most fundamental elements in power system operation, where its main objective is to ensure that the system is secure with respect to contingencies.

An optimal power flow (OPF) problem finds an optimal operating point for an objective function under given constraints. 
Security-constrained optimal power flow (SCOPF) concerns  contingency constraints within an OPF problem \cite{wood2013power}. It considers contingencies involving the disruption of generating stations and transmission lines, where controlling actions are carried out during power system operation and planning stage. Steady-state security is the capability of the system to function continuously within the rating of the equipment after a contingency has occurred \cite{alsac1974optimal}. The SCOPF problem is mainly classified into two classes: corrective and preventive formulations \cite{alsac1974optimal}, \cite{monticelli1987security}. The corrective formulation reschedules the power flow after the outage has occurred. The corrective SCOPF includes one set of control variables for normal operating conditions and another set of alteration variables for the contingency scenarios. 
The control variables in post-contingency are allowed to remove the violations caused by contingencies. The corrective model requires additional constraints and variables and many reschedules are required to perform the corrective actions. 
An iterative approach to solving corrective SCOPF is proposed in \cite{capitanescu2008new}. It considers a subset of potential contingencies, a steady-state security analysis, a contingency filtering technique, and checking post-contingency state feasibility \color{black}but it does not consider the parallel computation framework to speed-up computations. \color{black} On the other hand, preventive formulation considers normal situation control variables to minimize the cost function. The preventive SCOPF deals with one set of variables to fulfill both normal and contingency circumstances. A dispatch involving security constraints implements preventive control and thus better system security is achieved \cite{monticelli1987security}. The preventive generation rescheduling scheme using trajectory sensitivity analysis is proposed in \cite{nguyen2003dynamic} \color{black} but the algorithm is not robust. \color{black} 
An effective AC corrective/preventive contingency dispatch for the security-constrained unit commitment model to minimize the power system operational cost while preserving the system security is presented in \cite{fu2006ac} \color{black} but it does not consider parallel calculations for SCOPF subsystems. \color{black} The preventive SCOPF should only be applied when corrective actions are not enough or cannot be applied quickly \cite{wu2018robust}.

The main drawback of traditional SCOPF algorithms is the negligence of stability requirements and problem size because it is challenging to incorporate these in a forthright way \cite{capitanescu2011state}. Monte Carlo method is used for post-fault stability assessment and high-density sampling of operating points through time-domain simulations in \cite{konstantelos2016implementation}. 
A generation rescheduling method to increase the dynamic security of power systems is proposed in \cite{kuo1995generation}. \color{black} However, it does not model optimization problem and lacks various facets including prefault and postfault thermal and voltage limits of system operation. \color{black} 
A semi-definite programming relaxation to develop minimum singular value into voltage stability constrained OPF model is studied in \cite{wang2018sdp}. \color{black} A multi-agent predictive control model for load frequency control to increase the damping of oscillations in a power system using the Bat Inspired Algorithm is presented in \cite{elsisi2015multi}. The security of the smart meter data is analyzed using the machine learning techniques in \cite{elsisi2021reliable}. \color{black}
The recognition of critical control violations is very essential because their location in the system is critical for deciding on the installation of additional controls during the planning and the activation of the same controls during operation \cite{gunda2018remedial}. The solution of an optimization problem over a relaxed set gives lower bound on the optimal generation cost. However, the solution may risk system security and may not be feasible \cite{cui2018new}. A sparse tableau formulation for node-breaker representations in security-constrained optimal power flow is proposed in \cite{park2018sparse}. \color{black}However, it does not consider contingency filtering for large scale power systems and it leads to high computation burden.\color{black}

AGC is a speed governing property of a generation unit. It utilizes the power generation capacity to optimally distribute the power during disturbances in an interconnected system while minimizing the real-time generation cost. It offers significant cost savings under load variability and uncertain conditions while constraining the thermal rating limits of tie-lines. The AGC implemented in \cite{li2015connecting} ignores the tie-line thermal limits. The thermal limits on the tie-line are considered in \cite{mallada2017optimal} but AGC is not considered here. Besides, the corrective SCOPF presented in \cite{vrakopoulou2013probabilistic} does not consider the AGC response of generation units.

A large power system network involving numerous contingencies results in greater execution time and memory requirements. Thus, algorithms are created to consider the only potentially binding contingencies into the methodology \cite{capitanescu2007contingency}. 
The number of contingencies can be reduced using screening procedures \cite{ejebe1979automatic}. The network dimension reduction technique by ignoring the less affected areas of the system is proposed in \cite{karbalaei2018new}. Two contingency filtering techniques depending upon the post-contingency violations are proposed in \cite{capitanescu2007contingency}. The ranking procedure chooses contingencies based on a severity index that has a pre-defined threshold limit. 
A bi-level maximum-minimum optimization model to find the critical contingencies is developed in \cite{wu2018robust}. A contingency selection procedure based on the contingency explicit ranking is discussed in \cite{obadina1989var}. 


The security-constrained optimal power flow problem is a large scale optimization problem. One technique to solve SCOPF is Benders decomposition, which can take advantage of problem structure \color{black}\cite{li2008decomposed}\color{black}. Benders decomposition reduces the complexity of the problem by decomposing the original problem into master and various sub-problems \cite{shahidehopour2005benders}. 
Benders decomposition is used to solve the optimal power flow master problem in \cite{martinez2006security}. But the number of decision variables is greater for large power systems. It assumes short-term emergency ratings quite high that could result in voltage collapse or cascading outages before the corrective actions are effective. Benders decomposition is used to decompose the SCOPF problem into sub-problems associated with each contingency \cite{li2008decomposed}. 

There are three main challenges regarding the security-constrained ACOPF problem: considering full AC power flow constraints, modeling the AGC response, and scalability to tackle large scale problems with a large set of contingencies. For example, the contingency ranking method proposed in \cite{ni2003online} is not necessarily computationally efficient for a large network with numerous potential contingencies. Another example is the utilization of Benders decomposition to enhance the scalability  in \cite{street2010contingency} that considered DC power flow constraints, and not  AC power flow constraints.  

\subsection{Contributions}
\color{black} The main contributions of this paper are listed as follows:
\begin{enumerate}
    \item Modeling the AGC response of generations units in each contingency by introducing a set of valid constraints presenting the AGC response to changes in the system and integrate it with full AC power flow constraints. Thus, the procured solution presented in the security check sub-problem is the physical response of the system in case of contingency. 
    \item Integrating the convexified AC-OPF problem formulation in the Bender's decomposition  for both master problem and each contingency check sub-problems. The second-order cone relaxation of the master problem is employed to present the master problem in convex form. A tightness measure method is utilized to verify the merit of the employed convex relaxation approach. 
    \item Applying parallel processing techniques along with leveraging various level of contingency filtering to selected security check sub-problems. 
   The presented results verifies that compared to the traditional approaches, this presents significant improvement in solving problems with a large number of contingencies.  
\end{enumerate}\color{black}

\section{Problem Formulation} \label{problem_form}
 The original problem of SCOPF is presented in two parts which include the base problem and the constraints associated with each contingency. The base case problem formulation is given in (\ref{master_pblm}). The piece-wise objective function is represented by (\ref{mpa}), where h1 and h2 are linear function of bus injection mismatches, $c_g$ is the cost of generation units, and $\sigma_{ik}^{P\pm, Q\pm}$ are the mismatches due to contingencies on each bus. Base case bounds on voltage, real power, and reactive power are constrained by (\ref{mpb}), (\ref{mpc}), and (\ref{mpd}) respectively. Base case real and reactive power flows into a transmission line at the origin buses are defined by (\ref{mpe}) and (\ref{mpf}) respectively. The real and reactive power flow into a transformer at the origin buses in the base case is represented by (\ref{mpg}) and (\ref{mph}) respectively. Bus real and reactive power balance constraints in the base case are defined by (\ref{mpi}) and (\ref{mpj}) respectively. The base case line current rating at the origin bus is represented by (\ref{mpk}). The power ratings for the transformer in the base case at the origin bus are shown in (\ref{mpl}).

\begin{subequations}\label{master_pblm}
\begin{alignat}{1}
&\textbf{min} \sum\limits_{g\in \mathcal{G}} c_g+\frac{1}{\mid K\mid}\sum\limits_{k\in \mathcal{K}}(\sum\limits_{i}g(\sigma_{ik}^{P+}, \sigma_{ik}^{P-}, \sigma_{ik}^{Q+}, \sigma_{ik}^{Q-})+\sum\limits_{e\in \mathcal{E}}h_1(\sigma_{ek})+\sum\limits_{f\in \mathcal{F}}h_2(\sigma_{fk}))\label{mpa}\\
& \textbf{subject to:} \nonumber\\
& \underline{u}_i \leq u_i \leq  \overline{u}_i    \hspace{1cm}   \forall i \in \mathcal{I} \label{mpb}\\
& \underline{p}_g \leq p_g \leq  \overline{p}_g    \hspace{1cm}   \forall g \in \mathcal{G} \label{mpc}\\
& \underline{q}_g \leq q_g \leq  \overline{q}_g    \hspace{1cm}   \forall g \in\mathcal{G} \label{mpd}\\
& p_e^{ij} = g_{e}(e_{i}^2+f_{i}^2)-g_{e}(e_ie_j+f_if_j) -b_{e}(e_if_j-e_jf_i)\hspace{0.1cm}  \forall e \in \mathcal{E}\label{mpe}\\
& q_e^{ij} = -(b_e+b_e^{CH}/2)(e_{i}^2+f_{i}^2) +b_e (e_{i}e_{j}+f_{i}f_{j}) -g_e (e_if_j-e_jf_i) \hspace{0.1cm}   \forall e \in \mathcal{E}\label{mpf}\\
& p_{f}^{ij} = (g_f/\tau_f^2+g_f^M)({e_{i}}^2 + {f_{i}}^2)-g_f/\tau_f ({e_{i}}  {e_{j}} + {f_{i}}  f_{j})
-b_f/\tau_f ({e_{i}}  {f_{j}} - {e_{j}}  f_{i}) \hspace{0.2cm}   \forall f \in \mathcal{F}\label{mpg}\\
& q_{f}^{ij} = (b_f/\tau_f^2+b_f^M)({e_{i}}^2 + {f_{i}}^2)+b_f/\tau_f ({e_{i}}  {e_{j}} + {f_{i}}  f_{j})
-g_f/\tau_f ({e_{i}}  {f_{j}} - {e_{j}}  f_{i}) \hspace{0.2cm}   \forall f \in \mathcal{F}\label{mph}\\
& \sum\limits_{g\in G_i} p_g-p_i^L-g_i^{FS} u_i^2-\sum\limits_{j\in \delta(i)}p_e^{ij}
 -\sum\limits_{j\in \psi(i)}p_f^{ij}=0  \hspace{0.2cm}\forall i \in \mathcal{I} \label{mpi}\\
& \sum\limits_{g\in G_i} q_g-q_i^L-b_i^{FS}u_i^2-\sum\limits_{j\in \delta(i)}q_e^{ij}
 -\sum\limits_{j\in \psi(i)}q_f^{ij}=0  \hspace{0.2cm} \forall i \in \mathcal{I} \label{mpj}\\
& \sqrt{(p_e^{ij})^2+(q_e^{ij})^2} \leq \overline{R}_e  \hspace{0.4cm} \forall e \in \mathcal{E} \label{mpk}\\
&\sqrt{(p_f^{ij})^2+(q_f^{ij})^2} \leq \overline{s}_f \hspace{0.4cm} \forall f \in \mathcal{F} \label{mpl}
\end{alignat}
\end{subequations}
\vspace{-0.8cm}

The constraints associated with each contingency are given in (\ref{sub_pblm}). For each contingency, bounds on voltage, real power and reactive power are represented by (\ref{spa}), (\ref{spb}), and (\ref{spd}) respectively. The constraints (\ref{spc}) and (\ref{spe}) enforce the real and reactive power  of generators that are not active during contingency to zero.
The real and reactive power flow into a transmission line at the origin bus during a contingency are defined by (\ref{spf}) and (\ref{spg}) respectively. The real and reactive power flow into a transformer at the origin buses in each contingency are shown by (\ref{sph}) and (\ref{spi}) respectively. Bus real and reactive power balance constraints during a contingency along with soft constraint violation variables are represented by (\ref{spj}) and (\ref{spl}) respectively. The constraints (\ref{spk}) and (\ref{spm}) define that soft-constraint violation variables are positive. During each contingency, the line current rating at the origin bus with violation variable (\(\sigma_{ek}\)) is represented by (\ref{spn}) and (\ref{spo}). The transformer power ratings during a contingency at the origin bus along with constraint violation variables are shown in (\ref{spp}) and (\ref{spq}). The constraint (\ref{spr}) shows that an online generator but not selected to respond to contingency, retains its real output power from base case.
\vspace{-0.3cm}
\begin{subequations}\label{sub_pblm}
\begin{alignat}{2}
& \underline{u}_i^K \leq u_{ik} \leq \overline{u}_i^K \hspace{0.8cm} \forall k \in \mathcal{K}, i \in \mathcal{I} \label{spa}\\
& \underline{p}_g \leq p_{gk} \leq \overline{p}_g \hspace{0.8cm} \forall k \in \mathcal{K}, g \in \mathcal{G}_k \label{spb}\\
& p_{gk}=0 \hspace{1cm} \forall k \in \mathcal{K}, g \in \mathcal{G} \setminus \mathcal{G}_k \label{spc}\\
& \underline{q}_g \leq q_{gk} \leq \overline{q}_g \hspace{0.8cm} \forall k \in \mathcal{K}, q \in \mathcal{G}_k \label{spd}\\
& q_{gk}=0  \hspace{1cm} \forall k \in \mathcal{K}, g \in \mathcal{G} \setminus \mathcal{G}_k \label{spe}\\
& p_{ek}^{ij} = g_{e} ({e_{ik}^2}+{f_{ik}^2})-g_{e}(e_{ik}e_{jk}+f_{ik}f_{jk}) -b_{e}(e_{ik}f_{jk}-e_{jk}f_{ik})\hspace{0.1cm}   \forall k \in \mathcal{K}, e \in \mathcal{E}\label{spf}\\
& q_{ek}^{ij} = -(b_e+b_{e}^{CH}/2) ({e_{ik}^2}+{f_{ik}^2}) +b_e (e_{ik}e_{jk}+f_{ik}f_{jk})
-g_e (e_{ik}f_{jk}-e_{jk}f_{ik}) \hspace{0.8cm}   \forall k \in \mathcal{K}, e\in \mathcal{E}_k\label{spg}\\
& p_{fk}^{ij} = (g_f/\tau_f^2+g_f^M)({e_{ik}}^2 + {f_{ik}}^2)-g_f/\tau_f ({e_{ik}}  {e_{jk}} + {f_{ik}}  f_{jk})-b_f/\tau_f ({e_{ik}}  {f_{jk}} - {e_{jk}}  f_{ik}) \hspace{0.2cm}   \forall k \in \mathcal{K}, f \in \mathcal{F}_k\label{sph}
 \end{alignat}
\begin{alignat}{2}
& q_{fk}^{ij} = (b_f/\tau_f^2+b_f^M)({e_{ik}}^2 + {f_{ik}}^2)+b_f/\tau_f ({e_{ik}}  {e_{jk}} + {f_{ik}}  f_{jk})-g_f/\tau_f ({e_{ik}}  {f_{jk}} - {e_{jk}}  f_{ik}) \hspace{0.8cm}   \forall k \in \mathcal{K}, f \in \mathcal{F}_k\label{spi}\\
& \sum\limits_{g\in G_{ik}} p_{gk}-p_i^L-g_i^{FS} u_{ik}^2-\sum\limits_{j\in \delta_{k}(i)}p_{ek}^{ij}  -\sum\limits_{j\in \psi_{k}(i)}p_{fk}^{ij}=\sigma_{ik}^{P+}-\sigma_{ik}^{P-}  \hspace{0.8cm} \forall k \in \mathcal{K}, i \in \mathcal{I} \label{spj}\\
&\sigma_{ik}^{\rho} \geq 0  \hspace{1cm}   \forall k \in \mathcal{K}, i \in \mathcal{I}, ~\rho=\{\mathcal{P}\pm\} \label{spk}\\
& \sum\limits_{g\in G_{ik}} q_{gk}-q_i^L+b_i^{FS} u_{ik}^2-\sum\limits_{j\in \delta_{k}(i)}q_{ek}^{ij}
 -\sum\limits_{j\in \psi_{k}(i)}q_{fk}^{ij}=\sigma_{ik}^{Q+}-\sigma_{ik}^{Q-}  \hspace{0.8cm} \forall k \in \mathcal{K}, i \in \mathcal{I} \label{spl}\\
&\sigma_{ik}^{\rho} \geq 0  \hspace{1cm}   \forall k \in \mathcal{K}, i \in \mathcal{I}, ~\rho=\{Q\pm\} \label{spm}\\
& \sqrt{(p_{ek}^{ij})^2+(q_{ek}^{ij})^2} \leq \overline{R}_e^K+\sigma_{ek}^S \hspace{0.8cm} \forall k \in \mathcal{K},e \in \mathcal{E}_k \label{spn}\\
& \sigma_{ek}^{S} \geq 0   \hspace{1cm}   \forall k \in\mathcal{K}, e \in \mathcal{E}_k \label{spo}\\
& \sqrt{(p_{fk}^{ij})^2+(q_{fk}^{ij})^2} \leq \overline{s}_f^K+\sigma_{fk}^S \hspace{0.2cm} \forall k \in \mathcal{K}, f \in \mathcal{F}_k \label{spp}\\
& \sigma_{fk}^{S} \geq 0   \hspace{1cm}   \forall k \in \mathcal{K}, f \in \mathcal{F}_k \label{spq}\\
& p_{gk}=p_{g}  \hspace{1cm} \forall k \in \mathcal{K}, g \in \mathcal{G}_k \setminus \mathcal{G}_{^P}  \label{spr}
\end{alignat}
\end{subequations}
\vspace{-0.8cm}

The constraints presented in (\ref{sub_pblm}) do not cover the real and reactive power dispatch of generators participating in the contingency response. The real power (\(p_{gk}\)) of generator in a contingency $k$ is governed by the constraints given in (\ref{active_power_const}). 
Here, each generator has a pre-defined response rate given by $\alpha_g$ which is a portion of the total response given in $\Delta_k$ for the contingency $k$. This set of constraints represents the AGC response of the active generation units in the contingency, where the dispatch of generation units during the contingency $k$ is determined. If the determined value is outside of the physical limits of a generation unit, the dispatch is set to the physical limit and the AGC response is overwritten. The reactive power is subject to the constraints given in (\ref{reactive_power_const}), 
the voltage of buses with a generation bus is preserved to the values before the contingency if the reactive power limits are not reached.  
\vspace{-0.4cm}

\begin{equation}\label{active_power_const}
\begin{rcases}
& \{\underline{p}_{g}\leq p_{gk}\leq \overline{p}_{g} \hspace{0.2cm}\&\hspace{0.2cm}
p_{gk}=p_{g}+\alpha_g\Delta_k \}\\ 
& \{p_{gk}= \overline{p}_{g} \hspace{0.2cm}\&\hspace{0.2cm} p_{gk}\leq p_{g}+\alpha_g\Delta_k \} \\
& \{p_{gk}= \underline{p}_{g} \hspace{0.2cm}\&\hspace{0.2cm} p_{gk}\geq p_{g}+\alpha_g\Delta_k \}
\end{rcases}
k \in \mathcal{K}, g \in \mathcal{G}_k^P
\end{equation}
\vspace{-0.6cm}
\begin{equation}
\begin{rcases}
  & \{\underline{q}_{g}\leq q_{gk}\leq \overline{q}_{g} \&
  \sqrt{e_{{i_g}k}^2+f_{{i_g}k}^2}=\sqrt{e_{{i_g}}^2+f_{{i_g}}^2}\}\\
  & \{q_{gk}= \overline{q}_{g} \&\sqrt{e_{{i_g}k}^2+f_{{i_g}k}^2}\leq \sqrt{e_{{i_g}}^2+f_{{i_g}}^2} \}\\
  & \{q_{gk}= \underline{q}_{g} \&\sqrt{e_{{i_g}k}^2+f_{{i_g}k}^2}\geq \sqrt{e_{{i_g}}^2+f_{{i_g}}^2} \}
\end{rcases}k \in \mathcal{K}, g \in \mathcal{G}_k \label{reactive_power_const}
\end{equation}

These constraints represent the voltage control of the generation bus. If the upper limits of the reactive power generation are reached, the voltage on the connected bus is equal lower than that before the contingency. If the lower limits of the reactive power generation are reached, the voltage on the connected bus is equal greater than that before the contingency. The presented constraints in (\ref{active_power_const}) and (\ref{reactive_power_const}) are not straightforward to model in an optimization problem and will generally require employment of binary variables which cause an increase in the computation burden of each sub-problem. A set of valid constraints are introduced in the next section to solve the problem \emph{without employing binary variables which will adversely impact the computational efficiency of solving the problem}. 
\vspace{-0.5cm}

\section{Solution Method}\label{solution_meth}
\vspace{-0.1cm}
Solving the problem presented in section \ref{problem_form} has three challenges. First, the problem is nonlinear and non-convex. Thus, it is challenging to solve in polynomial time. Therefore, a convex relaxation technique is presented in subsection \ref{sssec:num1} to address this challenge. Second, the presented set of constraints in (\ref{active_power_const}) and (\ref{reactive_power_const}) are non-linear and cannot be directly incorporated into the AC-OPF problem formulation. Thus, a set of valid constraints are introduced in subsection \ref{agc_formulation} to address this challenge. Third, the set of contingencies dramatically increases the size of the problem. To address this challenge, a decomposition method tailored for the convexified problem is employed and it is discussed in subsection \ref{Benders_decomposition}.
 \vspace{-0.3cm}
\subsection{Formulating the Base Case Problem as a Second-Order Convex Relaxation} \label{sssec:num1}
 \vspace{-0.1cm}
The base case problem presented in (\ref{master_pblm}) is a non-convex optimization problem. The non-convexity arises because of the expressions \(e_{i}^2+f_{i}^2\),\hspace{0.1cm}\(e_ie_j+f_if_j\), and\hspace{0.1cm}\(e_if_j-e_jf_i\). To convexify the problem, new variables are defined as \(c_{ii}:=e_{i}^2+f_{i}^2\),\hspace{0.1cm}\(c_{ij}:=e_ie_j+f_if_j\) and\hspace{0.1cm}\(s_{ij}:=e_if_j-e_jf_i\). The revised base problem formulation is given in (\ref{form_mp}). The objective is to minimize the sum of the cost of generation and weight on the base case. Base case real and reactive power flow into a transmission line at the origin buses from (\ref{mpe}) and (\ref{mpf}) are updated to (\ref{subeq2}) and (\ref{subeq3}) respectively. The real and reactive power flow into a transformer at the origin buses in the base case from (\ref{mpg}) and (\ref{mph}) are revised to (\ref{subeq4}) and (\ref{subeq5}) respectively. The complex voltage magnitude at bus $i$ is restricted by the constraint (\ref{subeq6}). The characteristics of $c_{ij}$ and $s_{ij}$ are given in (\ref{subeq7}) and these new variables satisfy the relationship in (\ref{subeq8}) which is a second-order cone constraint. The base case problem is now convex and can be solved using off-the-shelf solvers.
\begin{subequations}\label{form_mp}
\begin{alignat}{2}
&\textbf{min} \sum\limits_{g\in G} c_g + \delta c^\sigma   \label{subeq1}\\
& \textbf{subject to: }  \nonumber \\
& (\ref{mpb})-(\ref{mpd})\nonumber\\
& p_e^{ij} = g_{e}c_{ii}-g_{e}c_{ij} -b_{e}s_{ij}\hspace{0.6cm}   \forall e \in \mathcal{E}\label{subeq2}\\
& q_e^{ij} = -(b_e+b_e^{CH}/2)c_{ii} +b_e c_{ij}-g_e s_{ij} \hspace{0.6cm}   \forall e \in \mathcal{E}\label{subeq3}\\
& p_f^{ij} = (g_e+g_{ii})/{t_m}^2 c_{ii} + (-g_e  t_r+b_e  t_i)/{t_m}^2 c_{ij} + (-b_e  tr-g_e  t_i)/{t_m}^2 s_{ij} \hspace{0.6cm} \forall f \in \mathcal{F}, e \in \mathcal{E}\label{subeq4}\\
& q_f^{ij} = -(b_e+b_{ii})/t_m^2 c_{ii}- (-b_e  t_r-g_e  t_i)/{t_m}^2 c_{ij} + (-g_e  t_r+b_e  t_i)/{t_m}^2 s_{ij} \hspace{0.6cm} \forall f \in \mathcal{F}, e \in \mathcal{E}\label{subeq5}\\
& (\ref{mpi})-(\ref{mpl})\nonumber\\
& \underline{u}_i^2 \leq c_{ii} \leq \overline{u}_i^2 \hspace{1cm} i \in \mathcal{I}   \label{subeq6}\\
& c_{ij}=c_{ji},\hspace{0.2cm} s_{ij}=-s_{ji} \hspace{1cm} (i,j) \in \mathcal{E}   \label{subeq7}\\
& c_{ij}^2+s_{ij}^2 \leq c_{ii}c_{ij} \hspace{1cm} (i,j) \in \mathcal{E}  \label{subeq8}
\end{alignat}
\end{subequations}
\vspace{-0.9cm}

\color{black}
One approach to check the gap between the original problem and the proposed method is to find the difference between the objective values and divide by the original problem objective. This method is used to check the gap of the solution in \cite{coffrin2018powermodels}. In second-order convex relaxation formulation, the constraint \eqref{subeq8} is in conic relaxation form. An alternative approach to check the tightness of the procured relaxation form is the difference between $c_{ij}^2+s_{ij}^2$ and $c_{ii}c_{jj}$ as given in \cite{soofi2020socp}. In a radial network, the difference is zero and the relaxation is exact due to the angles that can be uniquely identified by an arbitrary voltage reference. Nevertheless, in practice, the difference exists due to numerical precision issues but it is close enough to render a good quality solution. The quality of the solution is determined by measuring how far the difference is from zero. The negative logarithmic of the difference is taken due to the small value. Thus, the mathematical representation of tightness measure  given in (\ref{111}) can be leveraged as a post-process to determine the gap of the procured solution of the relaxed problem. Here, $(T)$ is the tightness measure for each pair $i$, $j$ and $c_{ii}$, $c_{jj}$, $c_{ij}$, $s_{ij}$ associated with pair $i$, $j$. Therefore, this gap indicator is applied in the algorithm as a measure of exactness for the procured solution from the convex relaxation problem. The higher the procured values, the tighter the relaxation is to the original non-convex problem.
\vspace{-0.1cm}
\begin{equation}
T_{ij} = -log \vert c_{ij}^2+s_{ij}^2-c_{ii}c_{jj}\vert 
\label{111}
\end{equation}
\vspace{-1.0cm}
\color{black}
\subsection{Modeling AGC Response}\label{agc_formulation}
This section addresses the generator real power contingency (\ref{active_power_const}) and reactive power contingency (\ref{reactive_power_const}) responses. Recurrent changes to the output of a power generator are necessary because the generated power and load should be balanced during a contingency in the power system. A set of valid constraints are presented in (\ref{form_sp}) for modeling the AGC response during a contingency.
\begin{subequations}\label{form_sp}
\begin{alignat}{2}
& (p_{gk}-(p_g+\alpha_{g}\Delta_{k}))(p_{gk}-\overline{p}_g)\leq 0   \hspace{0.8cm} \forall k \in \mathcal{K}, g \in \mathcal{G}_k^P \label{updated_a}\\
& (p_{gk}-(p_g+\alpha_{g}\Delta_{k}))(\underline{p}_g-p_{gk})\geq 0  \hspace{0.8cm} \forall k \in \mathcal{K}, g \in \mathcal{G}_k^P \label{updated_b}
\end{alignat}
\begin{alignat}{2}
& (c_{i_g,i_g}-c_{{i_g}k,{i_g}k})(q_{gk}-\overline{q}_g)\geq 0  \hspace{0.8cm} \forall k \in \mathcal{K},g \in \mathcal{G}_k^P \label{updated_c}\\
& (c_{{i_g}k,{i_g}k}-c_{i_g,i_g})(\underline{q}_g-q_{gk})\geq 0  \hspace{0.8cm} \forall k \in \mathcal{K},  g \in \mathcal{G}_k^P \label{updated_d}
\end{alignat}
\end{subequations}
The generator real power contingency response presented in (\ref{active_power_const}) is represented by (\ref{updated_a}) and (\ref{updated_b}). If the real power generation dispatch is in an open interval within the generation limits, the dispatch is determined using the participation factor and system wide changes of each contingency i.e. LHS of constraints (\ref{updated_a}) and (\ref{updated_b}) are both zero. If the real power generation dispatch hits the physical limits, the dispatch is equal to the physical limit and is not based on the participation factor and base case power dispatch. If the upper limits reached, LHS of (\ref{updated_a}) is zero and the second term on the LHS of (\ref{updated_b}) is negative, so its first term should be smaller than or equal to zero i.e the actual dispatch being equal or smaller than the desired one determined by the participation factor. Similarly, if the lower limits reached, LHS of (\ref{updated_b}) is zero and the second term on the LHS of (\ref{updated_a}) is negative, so its first term should be greater than or equal to zero i.e the actual dispatch being equal or greater than the desired one is determined by the participation factor. 

Moreover, reactive power contingency response presented in (\ref{reactive_power_const}) is revised to (\ref{updated_c}) and (\ref{updated_d}). The voltage of buses with a generator connected to them is equal to the base case voltage during a contingency if the reactive power dispatch of that generator is within an open interval of its physical limits, where the LHS of (\ref{updated_c}) and (\ref{updated_d}) are zero.  If the upper limits of the reactive power dispatch is reached, LHS of (\ref{updated_c}) is zero and the second term on the LHS of (\ref{updated_d}) is negative, so its first term should be smaller than or equal to zero i.e the generation bus voltage during the contingency is equal or smaller than the that of the base case. Similarly, if the lower limits of the reactive power dispatch reached, LHS of (\ref{updated_d}) is zero and the second term on the LHS of (\ref{updated_c}) is negative, so its first term should be smaller than or equal to zero i.e the generation bus voltage during the contingency is equal or greater than the that of the base case.  Thus, the procured constraints represent a set of valid constraints that can be integrated into the SCOPF problem formulation without introducing any binary variables. It should be noted that adding the equivalent constraints for the AGC response for each sub-problem will render a nonlinear problem. However, the complexity introduced by these constraints is not troublesome compared to the nonlinear inherited from AC power constraints that are treated with convex relaxation. These constraints presented in \eqref{updated_a}-\eqref{updated_b} and \eqref{updated_c}-\eqref{updated_d} are pairs which at least one of them is zero all the time. Beside, not all of these constraints are active in the solution process for solver.  Employing these constraints will significantly reduce the computation time for the AGC response in the SCOPF calculations as demonstrated in the case studies. 
\vspace{-0.45cm}

\subsection{Problem Decomposition}\label{Benders_decomposition}
\begin{figure}[t!] 
    \centering
    \includegraphics[width=0.55\textwidth]{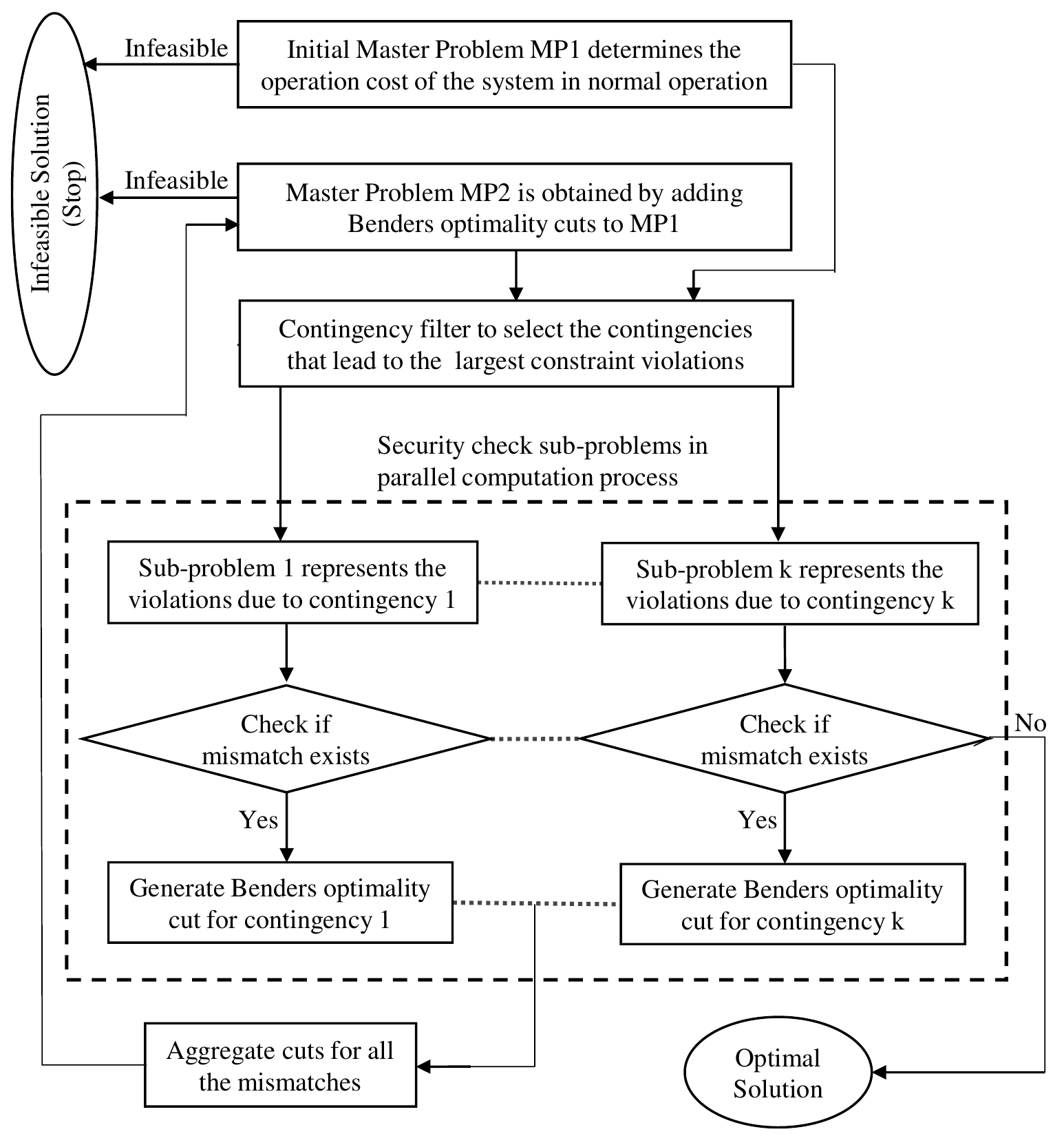}
    \caption{\color{black}Flowchart of Benders decomposition\color{black}}
    \label{block_dia} 
\end{figure}

To enable solving large-scale optimization problems, Benders decomposition is implemented and it is described in this subsection. 
The flowchart of the Benders decomposition is given in Fig. \ref{block_dia}. 

\begin{figure}[t!] 
    \centering
    \includegraphics[width=0.55\textwidth]{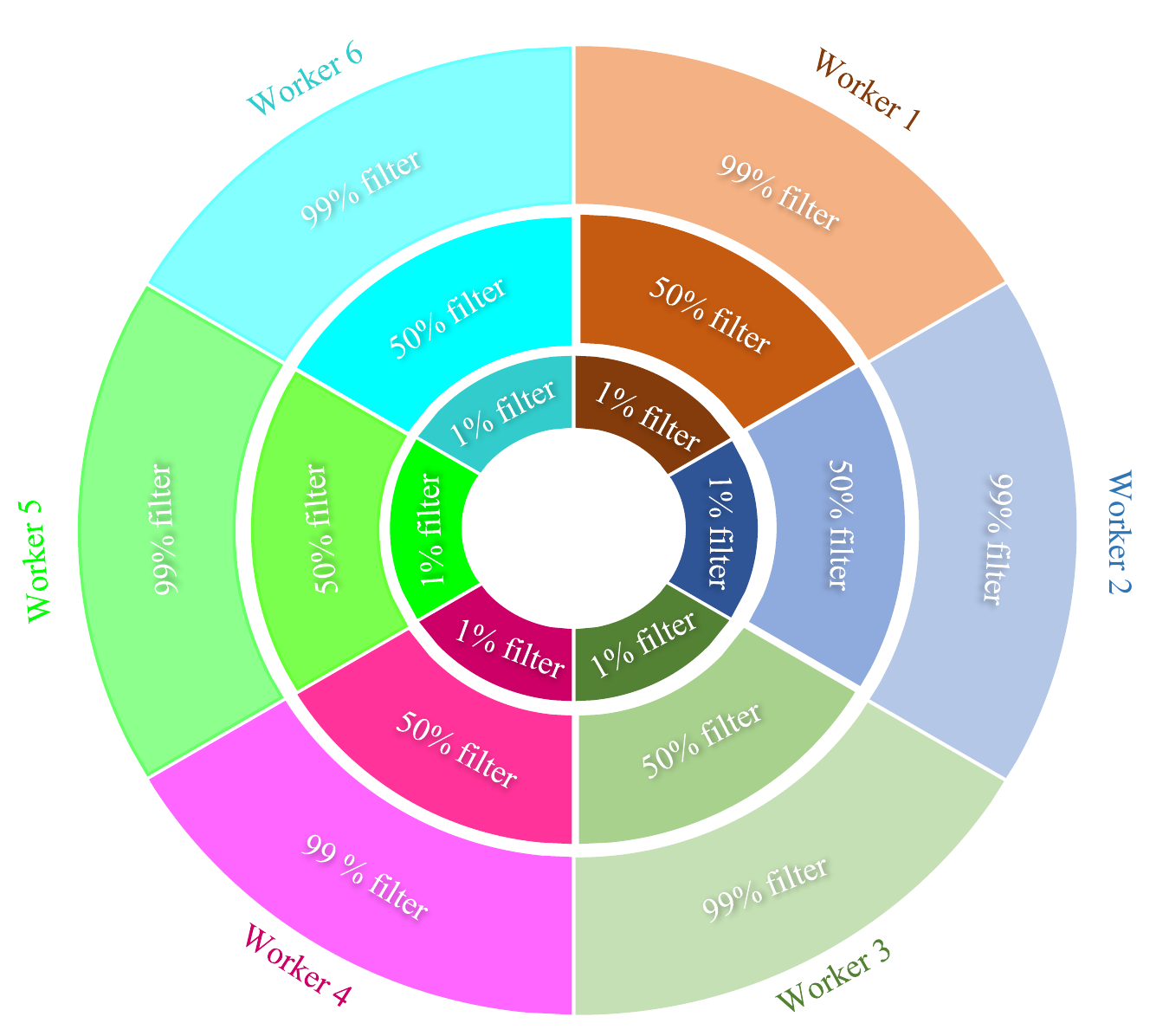}
    \caption{A sample donut chart of parallel computation process with 6 workers and three levels (1\%, 50\%, 99\%) of filtering.}
    \label{dia_pmap} 
\end{figure}

\begin{enumerate}
\item Solve the initial master problem MP1 using the convex relaxation method presented in subsection \ref{sssec:num1}. Obtain the lower bound (\(z_{lower}\)) for the objective value. If MP1 is infeasible, the original problem is infeasible. 

The abstract formulation of the master problem is given in (\ref{abs_mst_pblm}).
\vspace{-0.4cm}
\begin{subequations} \label{abs_mst_pblm}
\begin{alignat}{2}
&\underset{{x \in \mathcal{X}, u \in \mathcal{U} }}{min}\hspace{0.2cm}c^T x \label{subeq1}\\
& s.t. \hspace{0.5cm} A \textbf{x} +B \textbf{u}= d \hspace{1cm} \label{mp_equality_1}\\
& \hspace{0.5cm} \hspace{0.5cm} C \textbf{x} +D \textbf{u} \leq e \hspace{1cm}  \hspace{1cm}    \label{mp_inequality_1}
\end{alignat}
\end{subequations}
The equality constraint (\ref{mp_equality_1}) represents the equality constraints in (\ref{form_mp}) and the inequality constraint (\ref{mp_inequality_1}) represents the inequality constraints in (\ref{form_mp}).
\item Apply the contingency filtering technique to build a static ordering of all contingencies. The procured contingency filter ranks the components based on their significance in terms of utilization in the system. Thus, the contingencies that lead to the largest constraint violations are selected using the procured mismatches and a Benders cut for each violated contingency is added to the master problem. Three levels of filters are applied as shown in Fig. \ref{dia_pmap}, where 1\% filter represents the most critical contingencies. 

\item Pass the master problem solution to each sub-problem and solve it in a parallel computing process. The parallel computation technique accelerates the sub-problems solution because it involves multiple workers solving the security check sub-problems. It offers an environment that applies message passing interface to allow multiple workers to solve contingencies in a distributed memory. Each worker has its own private memory in a distributed memory environment. Message passing interface is a form of communication used in parallel computation that helps different workers to communicate with each other. In the proposed algorithm, there is one master worker and workers $1,...,n$, e.g. $n=6$  in Fig. \ref{dia_pmap}. The workers $1,...,n$ point their nodes towards the master worker. The master worker sends functions, data packets and signals to workers $1,...,n$ using message passing interface.  
The filtered contingencies are distributed among the available workers by dividing the total number of contingencies and the workers.  Finally, the workers return \(z_{upper}\) of the objective value for each sub-problem in the array.
The abstract composition of sub-problem is given in (\ref{abs_sub_pblm}).  
\vspace{-0.15cm}
\begin{subequations}\label{abs_sub_pblm}
\begin{alignat}{2}
& \underset{{x_k \in \mathcal{X}, u_k \in \mathcal{U} }}{min}\hspace{0.2cm}f(\bm{\sigma}) \\
& s.t. \hspace{0.5cm} F_k(\textbf{x}_k,\textbf{u}_k,{\sigma}_k)= 0 \hspace{1cm}  :\lambda_k  \label{sp_equality_1}\\
& \hspace{1cm} G_k(\textbf{x}_k,\textbf{u}_k) \geq 0 \hspace{1.3cm}  :\mu_k     \label{sp_inequality_1}
\end{alignat}
\end{subequations}
The subproblem objective is to minimize the mismatches. The equality constraint (\ref{sp_equality_1}) represents the equality constraints in (\ref{form_sp}) and the inequality constraint (\ref{sp_inequality_1}) represents the inequality constraints in (\ref{form_sp}).

The Benders cut for the master problem is given in (\ref{master_cut}).
\begin{equation}
 z_{lower}\geq d^{T}y+(h-Fy)^{T}\hat{u}^{p} \label{master_cut}
\end{equation}
\item Obtain the upper bound (\(z_{upper}\)) of the objective value from sub-problems solution. Check the mismatch between \(z_{upper}\) and \(z_{lower}\). If there is no mismatch, the optimal solution is procured. Otherwise, continue to the next step.
\item For each sub-problem, a Benders cut is generated. The cuts are aggregated and applied to the master problem MP2 with cuts. 
\item Repeat steps 3-5 until the total mismatch is less than a desired threshold.
\end{enumerate}

\vspace{-0.35cm}
\section{Case Study}\label{case_study}
In this section, three case studies are presented to demonstrate the effectiveness of the proposed algorithm for security-constrained AC-OPF. The first case study uses IEEE 14-bus system, the second uses 500-bus system, and the third uses 2000-bus system. The case studies are performed on a standard PC with an Intel Core i7-9700K CPU running at 3.60 GHz and 16.0 GB RAM. Julia built-in PMAP function is used for solving the sub-problems in parallel. PMAP applies the sub-problem constraints concurrently to the contingencies present in the array.
\vspace{-0.3cm}
\subsection{IEEE 14-Bus Power System}
In this case, a 14-bus power system, which consists of 14 buses, 5 generating units and 11 loads is deployed. This network contains 2 contingency scenarios where the outages occur on the branch connected between bus 6 and 12 and the least utilized generator 5 connected to bus 8. First, the master problem during the normal operating conditions is solved and the voltages at each bus are recorded as shown in column \(u_i\) of Table \ref{14_bus_voltages} corresponding to each bus index \(i\). For contingency scenario, the branch is taken out of system and the power flow is performed to identify the mismatches. The voltages during the branch contingency at each bus are shown by \(u_{i1}\). For the case of least utilized generator contingency, the voltages are shown in column \(u_{i2}\). The voltages of bus \(1\) for each contingency stays the same because it is a PV bus. This indicates the success of the presented solution method to preserve the voltage-controlled buses. The only exception is the voltage of bus \(8\) which is permitted to variate because the contingent generator is connected to it. The generator index column indicates the connection of each generator to the particular bus.

Table \ref{14_active_reactive} shows the automatic generation control response during each contingency. For the branch contingency, the real power generation \(p_{g1}\) from generators \(1-3\) and \(5\) is increased while it is the same for generator \(4\) because it was hitting the upper limit during the normal operating conditions. This indicates the success of the presented reformulation in preserving the AGC response of generation units under contingency. The reactive power \(q_{g1}\) for generators \(1, 3\) and \(5\) is increased while it is decreased for \(2\) and \(4\). As these reactive power values are within the limits of each generation unit, it is consistent with the voltage value presented in Table \ref{14_bus_voltages}. When the least utilized generator \(5\) is under contingency, the real and reactive powers from it are zero. In this case, the real power \(p_{g2}\) for generators \(1-3\) is increased while it is the same for generator \(4\)  because it reached the maximum limit. The reactive power \(q_{g2}\) is increased for generators \(1-3\) while it is decreased for generator \(4\). The minimum and maximum values of real and reactive power, as well as the participation factors for each generation unit, are shown in Table \ref{14_active_reactive}. The value of $\Delta_1$ is 0.02 MW for branch contingency and the value of $\Delta_2$ is -0.106 MW for generator contingency. \\

The tightness measure of the procured solution obtained from convex relaxation is determined using (\ref{111}) and tabulated in Table \ref{tightness_measure}. The gap between the proposed solution and the original problem is very small because the logarithmic value of the tightness is in-between 7 and 10. Thus, the tightness of the procured relaxation scheme represents that the proposed scheme is very close to the original one and a good quality solution is obtained.

\subsection{500-Bus Power System}
In this case, a 500-bus power system is utilized consisting of 90 generators, 200 loads and 131 transformers. It consists of 51 transformer contingencies and 326 branch contingencies. Table \ref{500_viol_conts} shows the results when the contingency filter is \(0.5\), which means that the top critical 50\% of total contingencies are taken into account. 
\color{black}The master problem objective value is written in \(0^{th}\) iteration and is determined without considering the violations in the optimization problem. This objective represents the operation cost $(\$)$ when the system operates in normal situation without any contingency. \color{black}With each iteration, the number of sub-problem contingencies that violated the rating decreases due to the Benders cuts applied to the master problem. An updated value of the objective is obtained after applying the cuts to the master problem. The objective value is increased from \(\$27,529\) to \(\$73,491\) in \(9\) iterations. The total mismatch for the contingencies decreases sharply with each iteration. The mismatch in \(2^{nd}\) iteration (21,379) is less than one half the value in \(1^{st}\) iteration (53,192). The algorithm converges in \(9^{th}\) iteration with \(0\) mismatch.  

By selecting critical contingencies active in the 500-bus power system, the effect on objective, the number of iterations, and the parallel (Par. Time) versus non-parallel (Series Time) solution convergence time of contingencies is tabulated in Table \ref{500_filter_change}. The percentage of active contingencies is decreased from \(99\%\) to \(N-1\). When \(99\%\) contingencies are active, the solution is converged in \(188\) iterations and convergence time is \(58,874\) seconds for parallel contingency solution case and very large for the non-parallel contingency solution case. For \(N-1\) contingency scenario, the solution is converged in \(3\) iterations and it took \(183\) seconds for the parallel and \(494\) seconds for the non-parallel situations. The merit of the parallel computation is that it solves the sub-problems on average of 18 times faster. For non-parallel cases, the computation time becomes infinitely large when large set of contingencies are evaluated. Thus, the contingency filtering technique ensures the security check of the procured AC-OPF solution by only considering the critical contingencies and the parallel computation approach make that happen in a time efficient way.

\subsection{2000-Bus Power System}
To demonstrate the efficiency of the proposed algorithm on a large power system network, a 2000-bus power system is considered. It consists of 544 generators, 1,125 loads, and 847 transformers. There are 432 generator contingencies and 2,753 branch contingencies summing to 3,185. In this case, top \(10\%\) of critical contingencies are considered. The mismatch, objective value, and contingencies that violated the power limits for each iteration are tabulated in Table \ref{2000_viol_conts}. 

The master problem objective value is \$741,696. During the first iteration, the objective value is increased to \$964,912 and the contingencies that violated the limits are 320. The value of mismatch is decreased sharply from \(352,940\) to \(0\) in \(5\) iterations. \color{black} It took 23 hours and 50 minutes for the algorithm to converge.  \color{black} The effectiveness of the proposed algorithm can be validated from the fact that after \(1^{st}\) iteration, the violated contingencies decreased from \(320\) to \(33\) and mismatch reduced from \(352,940\) to \(35,858\). \color{black}By considering 5\% critical contingencies, the algorithm converged in 12 hours and 18 minutes. \color{black}The algorithm is capable of solving critical contingencies from large number of contingencies in a time efficient way. Thus, the presented algorithm is computationally efficient for a large network with potential contingencies.

\vspace{-0.2cm}
\section{Conclusion}\label{conclusion}
This paper proposed a  reformulation for the automatic generation control in a decomposed convex relaxation algorithm to find the optimal solution to the ACOPF problem which is secure against a large number of contingencies. To solve this problem, the original ACOPF problem representing the system without contingency constraints is convexified by leveraging the second-order cone relaxation method. The contingencies are filtered to determine the corrective or preventive actions and the selected contingencies for  preventive action are considered in security constraints. Benders decomposition technique is utilized to decompose the security-constrained ACOPF into a master problem and several security check sub-problems. The sub-problems are evaluated in a parallel computing process to enhance the computational efficiency. AGC is modeled by a set of proposed valid constraints so that the solution obtained in each security check sub-problem is the physical response of each generation unit during contingency. Several case studies are presented to demonstrate the competence of the proposed valid AGC constraints and the scalability of the presented algorithm for the security-constrained ACOPF.
\color{black} An interesting future work is to convexify the automatic generation control response in the security-constrained AC optimal power flow calculations.\color{black} 
\vspace{-0.2cm}
\bibliographystyle{IEEEtran}


\bibliography{sample}

\begin{thebibliography}{10}
\providecommand{\url}[1]{#1}
\csname url@samestyle\endcsname
\providecommand{\newblock}{\relax}
\providecommand{\bibinfo}[2]{#2}
\providecommand{\BIBentrySTDinterwordspacing}{\spaceskip=0pt\relax}
\providecommand{\BIBentryALTinterwordstretchfactor}{4}
\providecommand{\BIBentryALTinterwordspacing}{\spaceskip=\fontdimen2\font plus
\BIBentryALTinterwordstretchfactor\fontdimen3\font minus
  \fontdimen4\font\relax}
\providecommand{\BIBforeignlanguage}[2]{{%
\expandafter\ifx\csname l@#1\endcsname\relax
\typeout{** WARNING: IEEEtran.bst: No hyphenation pattern has been}%
\typeout{** loaded for the language `#1'. Using the pattern for}%
\typeout{** the default language instead.}%
\else
\language=\csname l@#1\endcsname
\fi
#2}}
\providecommand{\BIBdecl}{\relax}
\BIBdecl

\bibitem{wood2013power}
A.~J. Wood, B.~F. Wollenberg, and G.~B. Shebl{\'e}, \emph{Power generation,
  operation, and control}.\hskip 1em plus 0.5em minus 0.4em\relax John Wiley \&
  Sons, 2013.

\bibitem{alsac1974optimal}
O.~Alsac and B.~Stott, ``Optimal load flow with steady-state security,''
  \emph{IEEE transactions on power apparatus and systems}, no.~3, pp. 745--751,
  1974.

\bibitem{monticelli1987security}
A.~Monticelli, M.~Pereira, and S.~Granville, ``Security-constrained optimal
  power flow with post-contingency corrective rescheduling,'' \emph{IEEE
  Transactions on Power Systems}, vol.~2, no.~1, pp. 175--180, 1987.

\bibitem{capitanescu2008new}
F.~Capitanescu and L.~Wehenkel, ``A new iterative approach to the corrective
  security-constrained optimal power flow problem,'' \emph{IEEE transactions on
  power systems}, vol.~23, no.~4, pp. 1533--1541, 2008.

\bibitem{nguyen2003dynamic}
T.~B. Nguyen and M.~Pai, ``Dynamic security-constrained rescheduling of power
  systems using trajectory sensitivities,'' \emph{IEEE Transactions on Power
  Systems}, vol.~18, no.~2, pp. 848--854, 2003.

\bibitem{fu2006ac}
Y.~Fu, M.~Shahidehpour, and Z.~Li, ``{A}{C} contingency dispatch based on
  security-constrained unit commitment,'' \emph{IEEE Transactions on Power
  Systems}, vol.~21, no.~2, pp. 897--908, 2006.

\bibitem{wu2018robust}
X.~Wu, A.~J. Conejo, and N.~Amjady, ``Robust security constrained
  {A}{C}{O}{P}{F} via conic programming: Identifying the worst contingencies,''
  \emph{IEEE Transactions on Power Systems}, vol.~33, no.~6, pp. 5884--5891,
  2018.

\bibitem{capitanescu2011state}
F.~Capitanescu, J.~M. Ramos, P.~Panciatici, D.~Kirschen, A.~M. Marcolini,
  L.~Platbrood, and L.~Wehenkel, ``State-of-the-art, challenges, and future
  trends in security constrained optimal power flow,'' \emph{Electric Power
  Systems Research}, vol.~81, no.~8, pp. 1731--1741, 2011.

\bibitem{konstantelos2016implementation}
I.~Konstantelos, G.~Jamgotchian, S.~H. Tindemans, P.~Duchesne, S.~Cole,
  C.~Merckx, G.~Strbac, and P.~Panciatici, ``Implementation of a massively
  parallel dynamic security assessment platform for large-scale grids,''
  \emph{IEEE Transactions on Smart Grid}, vol.~8, no.~3, pp. 1417--1426, 2016.

\bibitem{kuo1995generation}
D.-H. Kuo and A.~Bose, ``A generation rescheduling method to increase the
  dynamic security of power systems,'' \emph{IEEE Transactions on Power
  Systems}, vol.~10, no.~1, pp. 68--76, 1995.

\bibitem{wang2018sdp}
C.~Wang, B.~Cui, Z.~Wang, and C.~Gu, ``{S}{D}{P}-based optimal power flow with
  steady-state voltage stability constraints,'' \emph{IEEE Transactions on
  Smart Grid}, vol.~10, no.~4, pp. 4637--4647, 2018.

\bibitem{elsisi2015multi}
M.~Elsisi, M.~Aboelela, M.~Soliman, and W.~Mansour, ``Multi-agent model
  predictive control of nonlinear interconnected hydro-thermal system load
  frequency control based on bat inspired algorithm,'' \emph{optimization},
  vol.~1, no.~5, 2015.

\bibitem{elsisi2021reliable}
M.~Elsisi, K.~Mahmoud, M.~Lehtonen, and M.~M. Darwish, ``Reliable industry 4.0
  based on machine learning and iot for analyzing, monitoring, and securing
  smart meters,'' \emph{Sensors}, vol.~21, no.~2, p. 487, 2021.

\bibitem{gunda2018remedial}
J.~Gunda, G.~P. Harrison, and S.~Z. Djokic, ``Remedial actions for security
  constraint management of overstressed power systems,'' \emph{IEEE
  Transactions on Power Systems}, vol.~33, no.~5, pp. 5183--5193, 2018.

\bibitem{cui2018new}
B.~Cui and X.~A. Sun, ``A new voltage stability-constrained optimal power-flow
  model: Sufficient condition, {SOCP} representation, and relaxation,''
  \emph{IEEE Transactions on Power Systems}, vol.~33, no.~5, pp. 5092--5102,
  2018.

\bibitem{park2018sparse}
B.~Park, J.~Holzer, and C.~L. DeMarco, ``A sparse tableau formulation for
  node-breaker representations in security-constrained optimal power flow,''
  \emph{IEEE Transactions on Power Systems}, vol.~34, no.~1, pp. 637--647,
  2018.

\bibitem{li2015connecting}
N.~Li, C.~Zhao, and L.~Chen, ``Connecting automatic generation control and
  economic dispatch from an optimization view,'' \emph{IEEE Transactions on
  Control of Network Systems}, vol.~3, no.~3, pp. 254--264, 2015.

\bibitem{mallada2017optimal}
E.~Mallada, C.~Zhao, and S.~Low, ``Optimal load-side control for frequency
  regulation in smart grids,'' \emph{IEEE Transactions on Automatic Control},
  vol.~62, no.~12, pp. 6294--6309, 2017.

\bibitem{vrakopoulou2013probabilistic}
M.~Vrakopoulou, M.~Katsampani, K.~Margellos, J.~Lygeros, and G.~Andersson,
  ``Probabilistic security-constrained {A}{C} optimal power flow,'' in
  \emph{2013 IEEE Grenoble Conference}.\hskip 1em plus 0.5em minus 0.4em\relax
  IEEE, 2013, pp. 1--6.

\bibitem{capitanescu2007contingency}
F.~Capitanescu, M.~Glavic, D.~Ernst, and L.~Wehenkel, ``Contingency filtering
  techniques for preventive security-constrained optimal power flow,''
  \emph{IEEE Transactions on Power Systems}, vol.~22, no.~4, pp. 1690--1697,
  2007.

\bibitem{ejebe1979automatic}
G.~Ejebe and B.~F. Wollenberg, ``Automatic contingency selection,'' \emph{IEEE
  transactions on Power Apparatus and Systems}, no.~1, pp. 97--109, 1979.

\bibitem{karbalaei2018new}
F.~Karbalaei, H.~Shahbazi, and M.~Mahdavi, ``A new method for solving
  preventive security-constrained optimal power flow based on linear network
  compression,'' \emph{International Journal of Electrical Power \& Energy
  Systems}, vol.~96, pp. 23--29, 2018.

\bibitem{obadina1989var}
O.~Obadina and G.~Berg, ``{V}{A}{R} planning for power system security,''
  \emph{IEEE Transactions on Power Systems}, vol.~4, no.~2, pp. 677--686, 1989.

\bibitem{li2008decomposed}
Y.~Li and J.~D. McCalley, ``Decomposed {S}{C}{O}{P}{F} for improving
  efficiency,'' \emph{IEEE Transactions on Power Systems}, vol.~24, no.~1, pp.
  494--495, 2008.

\bibitem{shahidehopour2005benders}
M.~Shahidehopour and Y.~Fu, ``Benders decomposition: applying benders
  decomposition to power systems,'' \emph{IEEE Power and Energy Magazine},
  vol.~3, no.~2, pp. 20--21, 2005.

\bibitem{martinez2006security}
J.~Mart{\'\i}nez-Crespo, J.~Usaola, and J.~L. Fern{\'a}ndez,
  ``Security-constrained optimal generation scheduling in large-scale power
  systems,'' \emph{IEEE Transactions on Power Systems}, vol.~21, no.~1, pp.
  321--332, 2006.

\bibitem{ni2003online}
M.~Ni, J.~D. McCalley, V.~Vittal, and T.~Tayyib, ``Online risk-based security
  assessment,'' \emph{IEEE Transactions on Power Systems}, vol.~18, no.~1, pp.
  258--265, 2003.

\bibitem{street2010contingency}
A.~Street, F.~Oliveira, and J.~M. Arroyo, ``Contingency-constrained unit
  commitment with $ n-k $ security criterion: A robust optimization approach,''
  \emph{IEEE Transactions on Power Systems}, vol.~26, no.~3, pp. 1581--1590,
  2010.

\bibitem{coffrin2018powermodels}
C.~Coffrin, R.~Bent, K.~Sundar, Y.~Ng, and M.~Lubin, ``Powermodels. jl: An
  open-source framework for exploring power flow formulations,'' in \emph{2018
  Power Systems Computation Conference (PSCC)}.\hskip 1em plus 0.5em minus
  0.4em\relax IEEE, 2018, pp. 1--8.

\bibitem{soofi2020socp}
A.~F. Soofi, S.~D. Manshadi, G.~Liu, and R.~Dai, ``A socp relaxation for cycle
  constraints in the optimal power flow problem,'' \emph{IEEE Transactions on
  Smart Grid}, 2020.

\end{thebibliography}

\newpage

\clearpage

\clearpage

\begin{table}[h]
\setcounter{table}{0}
\small \centering
\caption{\color{black} Bus Voltages during Base and Contingency Cases } 
\begin{tabular}{ccccc}  \hline 
\textbf{Bus Index ($i$)}  & \bm{$u_i$}& \bm{$u_{i1}$} & \bm{$u_{i2}$} & \textbf{Generator Index} \\ \hline \hline
1 & \textbf{1.098} & \textbf{1.098} & \textbf{1.098}  & 1 \\\hline
2 & 1.093 & 1.093 & 1.093  & 2\\\hline
3 & 1.068 & 1.068  &1.068 & 3\\\hline
4 & 1.071 & 1.071 & 1.069  & -\\\hline
5 & 1.076 & 1.074 & 1.073  & - \\\hline
6 & 1.099 & 1.099 & 1.099  & 4\\\hline
7 & 1.095 & 1.089  & 1.080 & -\\\hline
8 & \textbf{1.099} & \textbf{1.099} & \textbf{1.080} & 5\\\hline
9 & 1.093 & 1.079 & 1.073  &  -\\\hline
10 & 1.087& 1.074  & 1.069  & - \\\hline
11 & 1.090 & 1.081 & 1.079  & -\\\hline
12 & 1.087 & 1.058 & 1.086  & -\\\hline
13 & 1.083 & 1.073  &1.079  & -\\\hline
14 & 1.073 & 1.060 & 1.059  & -\\\hline
\end{tabular}
\vspace{-0.25cm}
\label{14_bus_voltages}
\end{table} 

\clearpage

\begin{table}[h]
	\small \centering
	\caption{\color{black} Automatic Generation Control Response during a Contingency} 
	\begin{tabular}{cccccc} \hline 
		\textbf{Generator} \textbf{Index}  & \textbf{1} & \textbf{2} & \textbf{3} & \textbf{4} & \textbf{5}   \\ \hline \hline
		\bm{$p_g$}  
		& 37.96   & 83.29  & 5.80 & \textbf{110.50} & 0.32  \\\hline
		\bm{$q_{g}$}  
		& 1.84  & 18.69  & 27.62 & -4.98 & 2.85 \\\hline
		\bm{$p_{g1}$}   
	   & 38.07 & 83.67  & 6.79 & \textbf{110.50} &  0.38   \\\hline
		\bm{$q_{g1}$}  
	   & \textbf{2.57}  & 17.53 & 27.86 & -14.16 & \textbf{ 6.70}   \\\hline
		\bm{$p_{g2}$}   
		& 38.08 & 83.71  & 6.91 & \textbf{110.50} &  0   \\\hline
		\bm{$q_{g2}$}   
		& \textbf{3.16}  & 19.41  & 29.00 & -11.40 &  0   \\\hline
		\bm{$p_{min}$}  
		& 37.96  &  48.15 & 5.80 & 11.46 & 0.32    \\\hline
		\bm{$p_{max}$}
		& 245.45  &  157.50 & 82.50 & \textbf{110.50}& 80.14    \\\hline
		\bm{$q_{min}$}  
		& \textbf{-132.20}  & -76.98  & -0.88 & -19.11 & -8.03    \\\hline
		\bm{$q_{max}$}  
		& \textbf{76.06}  & 36.17  & 48.51 & 19.47 & 19.24     \\\hline
		\bm{$\alpha$}  
		& 5  & 19  & 49.3 & 38.8 & 3     \\\hline
	\end{tabular}
	\label{14_active_reactive}
\end{table}

\clearpage

\vspace{-0.4cm}
\begin{table}[h]
\small \centering
\caption{\color{black} Tightness Measure of Bus-pairs of the IEEE 14-bus System} 
\begin{tabular}{cccccc}  \hline 
\textbf{Bus-pair}  & \bm{$T_{ij}$}& \textbf{Bus-pair}  & \bm{$T_{ij}$}& \textbf{Bus-pair}  & \bm{$T_{ij}$}\\ \hline \hline
(1,2) & 8.0012 & (6,11) & 7.7436 &(5,6) & 8.1014\\\hline
(1,5) & 8.0060 & (6,12) & 7.7453& (13,14) & 8.0052\\\hline
(2,3) & 8.0053& (6,13) & 7.7453 &(4,9) & 7.7431\\\hline
(2,4) & 8.0036 & (7,8) & 7.6488& (12,13) & 8.0024\\\hline
(2,5) & 8.0036 & (7,9) & 8.7622 &(4,7) & 9.8836\\\hline
(3,4) & 8.0029 &(9,10) & 8.0015 &(10,11) & 8.0032\\\hline
(4,5) & 8.0008& (9,14) & 8.0042 &\\\hline
\end{tabular}
\label{tightness_measure}
\end{table}

\clearpage

\begin{table}[h]
	\small \centering
	\caption{\color{black} Summary of Results for 500-Bus System } 
	\begin{tabular}{cccc} \hline 
		\textbf{Iterations [n]}  & \textbf{Violations} & \textbf{Objective [\$]}  & \textbf{Mismatch} \\ \hline \hline
		0 
		&  -  & 27,529.7  & - \\\hline
		1 
		& 188   & 73,491.5  & 53,192.4 \\\hline
		2 
		& 93 & 73,491.5  & 21,379.8  \\\hline
		3  
		& 46  & 73,491.5  & 9,164.1  \\\hline
		4  
		& 22  & 73,491.5  & 4,812.1  \\\hline
		5  
		& 11  & 73,491.5  & 2,379.3  \\\hline
        6 
		& 5 & 73,491.5  & 937.8  \\\hline
		7  
		& 2  & 73,491.5  & 362.0  \\\hline
		8  
		& 1  & 73,491.5  & 170.5  \\\hline
		9  
		& 0  & 73,491.5  & 0  \\\hline
	\end{tabular}
	\label{500_viol_conts}
\end{table}

\clearpage

\begin{table}[h]
	\small \centering
	\caption{\color{black} Different Contingency Filters to observe the effect on Objective Value, Iterations and Solution Convergence Time} 
	\vspace{-0.2cm}
	\begin{tabular}{ccccc} \hline 
		\textbf{Filter} & \textbf{Obj. [\$]} & \textbf{Iters.} & \textbf{ Par. Time [s]} & \textbf{ Series Time [s]}  \\ \hline \hline
		0.01 
		& 131,234.2  & 188  & 58,874 &N/A\\\hline
		0.25 
		& 93,796.6 & 18  & 4,941 &  139,015\\\hline
		0.5  
		& 73,491.5  & 9  & 2,123 & 46,814\\\hline
		0.75  
		& 45,849.1  & 5  &  765 & 16,053\\\hline
		0.99  
		& 28,333.5  & 3 &  183 & 494\\\hline 
	\end{tabular}
	\label{500_filter_change}
\end{table}
\vspace{-0.2cm}

\clearpage

\hspace{-0.8cm}
\begin{table}[h]
	\centering
	\caption{\color{black} Summary of Results for 2000-Bus System } 
	\begin{tabular}{cccc} \hline 
		\textbf{Iterations [n]} & \textbf{Violations} & \textbf{Objective [\$]} & \textbf{Mismatch} \\ \hline \hline
		0 
		&  -  & 741,696.3  &  -\\\hline
		1 
		&  320  & 964,912.7  &  352,940.7\\\hline
		2 
		& 33 & 964,912.7  & 35,858.7  \\\hline
		3  
		&  3 & 964,912.7  &  3,291.1 \\\hline
		4  
		&  1 & 964,912.7  &  936.9 \\\hline
		5  
		&  0 & 964,912.7  &  0 \\\hline
	\end{tabular}
	\label{2000_viol_conts}
\end{table}

\end{document}